\documentclass{amsart}
\usepackage{amsmath,amsthm,amssymb,IMjournal}

\begin{document}
\newtheorem{The}{Theorem}[section]

\numberwithin{equation}{section}

\title{A note on the independent  domination number versus the domination number in  bipartite graphs}

\author{\|Shaohui |Wang|, Garden City,
        \|Bing |Wei|, University}

\rec {Acceptted by Czechoslovak Mathematical Journal}

\dedicatory{The first  author is
partially supported by the Summer Graduate Research Assistantship Program of Graduate School at the University of Mississippi.}

\abstract 
Let $\gamma(G)$ and $i(G)$ be the domination number and the independent domination number of $G$, respectively. 
Rad and Volkmann  posted a conjecture that $i(G)/ \gamma(G) \leq \Delta(G)/2$ for any graph $G$, where $\Delta(G)$ is its maximum degree (See  \cite{5}: N.J. Rad, L. Volkmann, A note on the independent domination number in graphs. Discrete Appl. Math. 161(2013) 3087--3089).  In this work, we verify the conjecture for  bipartite graphs. Several  graph classes attaining the extremal bound and graphs containing odd cycles with the ratio larger than $\Delta(G)/2$
are provided as well.
\endabstract



\thanks
  The first  author is
partially supported by the Summer Graduate Research Assistantship Program of Graduate School at the University of Mississippi.
\endthanks

Let $G= (V(G), E(G))$ be a simple  undirected graph on $n$ vertices.
    For  $v\in V(G)$,  $N_G(v)=\{ w\in V(G):  vw\in E(G)\}$ is 
an  open neighborhood of $v$ in $G$.  If $N_G(v) = \phi$,   $v$ is called an isolated vertex.    For $S\subseteq V(G)$, $N_G(S)$ is an open neighborhood of $S$ and  $G-S$ is a subgraph induced by $V(G)-S$.  
A vertex set $D \subseteq V(G)$  is a  dominating set if every vertex of $V(G)-D$ is adjacent to some vertices of $D$. The minimum cardinality of a dominating set of $G$ is called the domination number, denoted by $\gamma(G)$.  Furthermore,
a vertex set $I \subseteq V(G)$  is  an independent dominating set if $I$ is both an independent set and a dominating set in $G$, where an independent set  is a set of vertices in a graph such that  no two of which are adjacent. The minimum cardinality of  an independent dominating set of $G$ is called the  independent domination number, denoted by $ i(G)$. We refer to \cite{10, 11} for notation and terminologies used but not defined here. Currently, lots of  work related to domination number and independent domination number have been done, see  \cite{1978, 2, 1977}. 

In general, it is very difficult  to find the domination and independent domination numbers. Note that $i(G) \geq \gamma(G)$, which implies $i(G) / \gamma(G) \geq 1$. A natural question is to determine an upper
bound for $i(G)/ \gamma(G)$.
Hedetniemi and   Mitchell \cite {1977}   in 1977  showed that if $L$ is a line graph of a tree, then
 $ i(L) / \gamma (L) =1$,   where
the line graph $L(G)$ of  a connected graph $G$ is a graph such that each vertex of $L(G)$ represents an edge of $G$  and two vertices of $L(G)$ are adjacent if and only if their corresponding edges share a common endpoint in $G$.
Since a line graph does not have an induced subgraph isomorphic to $K_{1,3}$, 
 Allan and Laskar \cite {1978} extended the previous result  and obtained that if a graph $G$ is a  $K_{1,3}$-free graph, then $  i(G)/ \gamma(G) = 1$. 
In 2012, Goddard et al.\cite{3}  continued the similar approach and proved that $i(G)/ \gamma(G) \leq  3/2$ if $G$ is a cubic graph. In 2013, Southey and Henning \cite{6} improved the previous bound to $i(G)/ \gamma(G)  \leq 4/3$ for a connected cubic graph $G$ other than $K_{3,3}$.  
Additionally, Rad and Volkmann \cite{5}  obtained an upper bound of $i(G)/\gamma(G)$ related to the maximum degree $\Delta(G)$ for any graph $G$ and proposed a conjecture below.
 \begin{The}\label{T1} { (Rad and Volkmann \cite{5}) 
If $G$ is a graph, then
$$ \frac{i(G)}{\gamma(G)} \leq
 \left\{ \begin{array}{rcl}
   \frac{\Delta(G)}{2},\;\;\;\;\;\;\;\; \;\;\;\;\;\;\;\;\;\;\;\;\;\;\;&& \mbox{if } 3 \leq \Delta(G) \leq 5,
  \\\Delta(G) -3 +\frac{2}{\Delta(G) - 1},&&   \mbox{if } \Delta(G) \geq 6.
\end{array}\right.$$
} \end{The}
  \begin{The}\label{T2} \emph{\cite{5}  
If $G$  is a graph with $\Delta(G) \geq 3$, then $i(G)/\gamma(G) \leq \Delta(G)/2$.
}
\end{The}

In 2014, Furuya et al.\cite{7} proved that  $i(G)/ \gamma(G)  \leq \Delta(G) - 2 \sqrt{\Delta(G)} +2$ and gave a class of graphs which achieve the new upper bound.  However,
when $\Delta(G) \neq  4$,   $ \Delta(G) - 2 \sqrt{\Delta(G)} +2 > \Delta(G) / 2.$  On the other hand, it is still very interesting to determine other class of graphs, for which Conjecture 2 holds.

Motivated by Conjecture 2 and previous results, we will show that:

 \begin{The}\label{T3} \emph{
If $G$ is a bipartite graph with maximum degree $\Delta(G) \geq 2$, then $$ \frac{i(G)}{\gamma(G)} \leq
\frac{\Delta(G)}{2}.
$$ 
} \end{The}

We now provide graphs containing some odd cycles, for which Conjecture 2 does not hold. For any large integer $n$,
the graph $G'$ consists of an odd cycle $C_{2k+1}$ and $(2k+1)s$  vertices of degree $1$ such that each vertex on $C_{2k+1}$ is adjacent to exactly $s$ vertices of degree $1$, 
for any positive integers $k, s $.
Then $\Delta(G') = s+2, \gamma(G') = 2k+1$ and $i(G') = k+(k+1)s$. By  simple calculations, we can get $i(G') / \gamma(G') > \Delta(G') / 2$ if $s > 2k+2$.

 \vskip 2mm \noindent{\bf Proof of Theorem 3.}
Let $A$ and $B$ be the partite sets of the bipartite graph $G$, and  $D$ be a minimum dominating set of $G$. Assume that  $I_0$ is the set of isolated vertices in $G[D]$. Set $A_0 = A \cap I_0$, $A_1 = ( D \setminus A_0) \cap A$,  $B_0 = B \cap I_0$ and $B_1 = (D \setminus B_0) \cap B$.
Then $|D|= |A_0| + |A_1|+|B_0|+|B_1|$. Renaming the sets $A$ and $B$, if necessary, we may assume that $|A_1| \geq |B_1|$, which implies that $|B_1| \leq (|D|-|A_0|-|B_0|)/2 \leq |D|/2$. Let $A_2 = A \setminus (A_0 \cup A_1 \cup N_G(B_0))$. Since $D$ is a dominating set of $G$, we have that $A_2 \subset N_G(B_1)$. 
Furthermore, by the choice of $I_0$, any vertex of $A_1$(or $B_1$, respectively) is adjacent to at least one vertex of  $B_1$(or $A_1$, respectively) in $G[D]$. So, $|N_G(B_1) - A_1 - N_G(B_0)| \leq (\Delta -1)|B_1|$, which implies that $|A_2| \leq (\Delta -1)|B_1|$.
We now consider the set $I= A_0 \cup A_1 \cup A_2 \cup B_0$. By the constructions, we see that $I$ is an independent set of $G$. We now show that $I$ is a dominating set. In fact, if $v \in A \setminus I$, then $v$ is dominated by $B_0$, and  if $v \in B \setminus I$, then either $v \in B_1$, in which case $v$ is dominated by $A_1 \cup A_2$, or $v \in B \setminus (B_0 \cup B_1)$, in which case $v$ is dominated by $A_0 \cup A_1$. Thus, $I$ is an independent dominating set of $G$, and so
 \begin{eqnarray}
i(G) &\leq &  |I| \nonumber \\
&= &  |A_0|+|A_1|+|A_2|+|B_0| \nonumber \\
&=&  |D| -|B_1|+|A_2| \nonumber \\
&\leq &  |D|-|B_1|+(\Delta -1)|B_1| \nonumber \\
&= &  |D|+|B_1|(\Delta -2) \nonumber \\
&\leq &  |D| + \frac{1}{2}|D|(\Delta -2) \nonumber \\
&= &  \frac{1}{2}|D| \cdot \Delta \nonumber \\
&= &  \frac{1}{2} \gamma (G) \cdot \Delta, \nonumber
\end{eqnarray}

or,  equivalently, $\frac{i(G)}{\gamma(G)} \leq \frac{\Delta}{2}$.
 $\hfill\Box$
 \vskip 2mm \noindent
{\bf Remarks.}
We see that Conjecture 2 holds for the bipartite graph $G$. The upper bound  $\frac{\Delta}{2}$ can be achieved, and a balanced double star and a complete balanced bipartite graph are examples of bipartite graphs to attain the upper bound. 
\vskip 0.2cm
 \vskip 2mm \noindent{\bf Acknowledgements.}
The authors would like to thank the referees for their valuable comments, which lead to shorten the main proof  in this paper.

{\small
{\em Authors' addresses}:
{\em Shaohui Wang}, Department of Mathematics, The University of Mississippi, University, MS 38677, USA;  Department of Mathematics and Computer Science, Adelphi University, Garden City, NY 11530, USA,
 e-mail: \texttt{shaohuiwang@yahoo.com, swang4@go.olemiss.edu}.

{\em Bing Wei}, Department of Mathematics, The University of Mississippi, University, MS 38677, USA,
 e-mail: \texttt{bwei@olemiss.edu}.

}

\end{document}